\documentclass[reqno]{gokova}
\usepackage{epsfig,graphicx}

\begin{document}
\def\E{\ifmmode{\mathbb E}\else{$\mathbb E$}\fi} 
\def\N{\ifmmode{\mathbb N}\else{$\mathbb N$}\fi} 
\def\R{\ifmmode{\mathbb R}\else{$\mathbb R$}\fi} 
\def\Q{\ifmmode{\mathbb Q}\else{$\mathbb Q$}\fi} 
\def\C{\ifmmode{\mathbb C}\else{$\mathbb C$}\fi} 
\def\H{\ifmmode{\mathbb H}\else{$\mathbb H$}\fi} 
\def\Z{\ifmmode{\mathbb Z}\else{$\mathbb Z$}\fi} 
\def\P{\ifmmode{\mathbb P}\else{$\mathbb P$}\fi} 
\def\T{\ifmmode{\mathbb T}\else{$\mathbb T$}\fi} 
\def\SS{\ifmmode{\mathbb S}\else{$\mathbb S$}\fi} 
\def\DD{\ifmmode{\mathbb D}\else{$\mathbb D$}\fi} 

\renewcommand{\a}{\alpha}
\renewcommand{\b}{\beta}
\renewcommand{\d}{\delta}
\newcommand{\D}{\Delta}
\newcommand{\e}{\varepsilon}
\newcommand{\g}{\gamma}
\newcommand{\G}{\Gamma}
\newcommand{\la}{\lambda}
\newcommand{\La}{\Lambda}
\newcommand{\n}{\nabla}
\newcommand{\var}{\varphi}
\newcommand{\s}{\sigma}
\newcommand{\Sig}{\Sigma}
\renewcommand{\t}{\tau}
\renewcommand{\th}{\theta}
\renewcommand{\O}{\Omega}
\renewcommand{\o}{\omega}
\newcommand{\z}{\zeta}

\newcommand{\ben}{\begin{enumerate}}
\newcommand{\een}{\end{enumerate}}
\newcommand{\be}{\begin{equation}}
\newcommand{\ee}{\end{equation}}
\newcommand{\bea}{\begin{eqnarray}}
\newcommand{\eea}{\end{eqnarray}}
\newcommand{\bc}{\begin{center}}
\newcommand{\ec}{\end{center}}

\newcommand{\IR}{\mbox{I \hspace{-0.2cm}R}}
\newcommand{\IN}{\mbox{I \hspace{-0.2cm}N}}

\newtheorem{thm}{Theorem}[section]
\newtheorem{cor}[thm]{Corollary}
\newtheorem{lem}[thm]{Lemma}
\newtheorem{prop}[thm]{Proposition}
\newtheorem{ax}{Axiom}
\newtheorem{conj}[thm]{Conjecture}

\theoremstyle{definition}
\newtheorem{defn}{Definition}[section]

\theoremstyle{remark}
\newtheorem{rem}{\rm\bfseries{Remark}}[section]
\newtheorem*{notation}{Notation}

\newtheorem{ques}{\rm\bfseries{Question}}[section]
\newtheorem{cons}[rem]{\rm\bfseries{Construction}}
\newtheorem{exm}[rem]{\rm\bfseries{Example}}



\title{Variations on Fintushel-Stern Knot Surgery on $4$-manifolds }
\author[AKBULUT]{ Selman Akbulut $^1$}
\address{ MSU, Dept of Mathemtics, E.Lansing, MI, 48824, USA}
\email{akbulut@math.msu.edu}
\thanks{$^1$Supported in part by NSF fund DMS 9971440 }

\begin{abstract}
We discuss some consequences Fintushel-Stern `knot surgery' operation 
on $4$-manifolds coming from its handlebody description. We give
some generalizations of this operation and give a counterexample to their 
conjecture. 
 \end{abstract}

\volume{8}

\maketitle

\section{Introduction}

Let $X$ be a smooth $4$-manifold and $ K\subset S^{3} $ be a knot, In \cite{FS} among other
things Fintushel and Stern had shown that the operation $ K\to X_{K} $ of
replacing a tubular neighborhood of imbedded torus in $X$ by $ (S^{3} - K)\times
S^{1}$ could results change of smooth structure of $X$. In \cite{A} an algorithm of
describing handlebody of $X_{K}$ in terms of the handlebody of $X$ was described. In
this article we will give some corollaries of this construction, and present a
counterexample to conjecture of Fintushel and Stern which was overlooked in \cite{A}.
First we need to recall the precise description of $X_{K}$: Recall that the first
picture of Figure 1 is $ T{^2} \times D^{2} $, and the second one is the cusp $C$  
(i.e. $ B^{4}$ with a
$2$-handle attached along the trefoil knot with the zero framing). Clearly the
cusp $C$ contains a copy of $ T{^2} \times D^{2} $.

\begin{figure}[htb]
\includegraphics{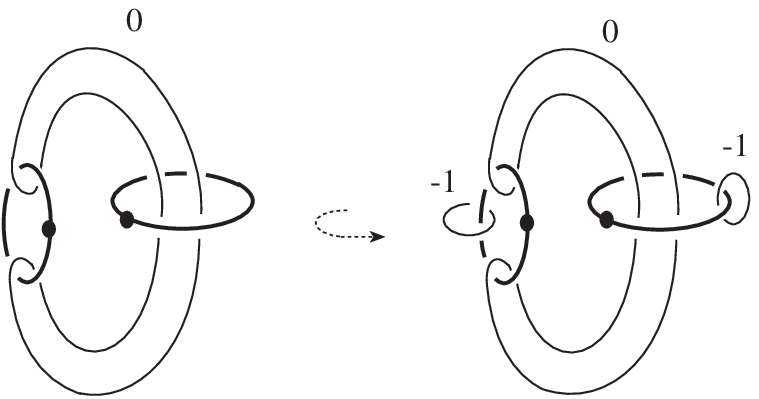}
\caption{}
\label{myfigure}
\end{figure}

 In \cite{FS} an imbedded torus $ T^{2}\subset X $ is called a {\it c-imbedded}
torus if it has a cusp neighborhood in $X$, i.e. $ T^{2}
\hookrightarrow C\hookrightarrow X $ as in Figure 1.  Now let  $ N \approx  K\times D^2 $ be
the trivialization of the open tubular neighborhood of the knot $K$ in $S^3$ given by
the
$0$-framing. Let 
$ \varphi: \partial (T^{2}\times D^{2})\to \partial  (K\times D^{2}) \times
S^{1}
$ be any diffeomorphism with $ \; \varphi (p\times \partial D^2)=K\times p $,
where
$ p\in T^{2} $ is a point, then define: 
$$X_{K}= (X-T^2\times D^2)
\smile_{\varphi} (S^{3}-N)\times S^1$$ 

Let $ Spin_{c}(X) $ be the set of $ Spin_{c} $ structures on $X$, e.g. if
$H_{1}(X)$ has no $2$-torsion then.
$$Spin_{c}(X) =\{ \; a \in H^{2}(X;{\bf Z}) \; |\; a=w_{2}(TX) \;\mbox{mod} 2 \;\} $$
Recall that Seiberg-Witten invariant $SW_{X}$ of $X$ is a
symmetric function 
$$SW_{X}:Spin_{c}(X)\to {\bf Z} $$
It is known that the function
$ SW_{X} $ is nonzero on the complement of a finite set ${\it B}=\{\pm \a_1, \pm
\a_2,..,\pm a_n\}$ which is called {\it the set of basic homology classes classes}. By
setting $ \a_{0} =0$ and $t_{j}=exp(\a_{j})$, the function $ SW_{X} $  is usually
written as a single polynomial
$$ SW_{X}= \sum_{j=0}^{n}SW_{X}(\a_{j})t_{j} $$
 
Now if $T$ is a $c$-imbedded torus in $X$ , and  $[T]$ be the homology class in
$H_{2}(X_{K};{\bf Z})$ induced from  $T^{2}\subset X$, and $ t=\mbox{exp}(2[T]) $, and
$\Delta_{K}(t)$ the Alexander polynomial of the knot $K$ (as a symmetric
Laurent polynomial), then Fintushel and Stern \cite{FS} theorem says:

\begin{thm}
$ SW_{X_{K}}=SW_{X} . \;\Delta_{K}(t)$
\end{thm}

\vspace{.15in}

Recall that in \cite{A} the algorithm of drawing the handlebody of $X_{K}$ from $X$ is
described  as follows: First we identify the {\it core circles} of the 1-handles of the
handlebody of $S^{3} - K$

\begin{figure}[htb]
\includegraphics{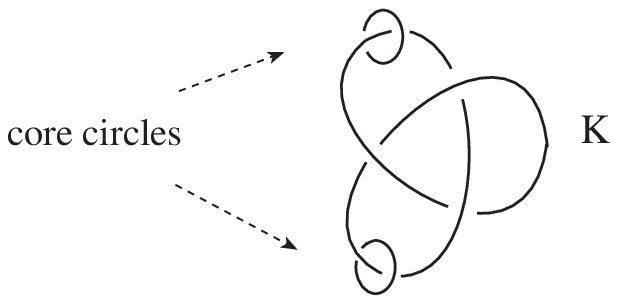}
\caption{}
\label{myfigure}
\end{figure}

Then when we see an imbedded cusp $ C $ in the handlebody of $X$  as in the first
picture of Figure 3, we change it to the second picture $\tilde{C}$ of Figure 3. 
This means
that we change one of the
$1$-handles of $T^{2}\times D^{2}$ inside of $C$ to the ``slice $1$-handle" 
obtained from $K\#  (-K)$ (i.e. remove the obvious slice disk which
$ K \#  (-K) $ bounds from $B^4$), and connect the {\it core circles} 
of the knots $K$ and $-K$ by $2$-handles as
shown in Figure 3. More precisely, there is a diffeomorphism 
between the boundaries of manifolds $C$ and
$\tilde{C}$ of Figure 3, and the operation $X\to X_{K}$ 
corresponds to cutting out $ T^{2}\times D^2 $ from
$X$ and gluing the second manifold of Figure 3 by this 
diffeomorphism (in the figure $K$ is drawn as the
trefoil not). 

\begin{figure}[htb]
\includegraphics{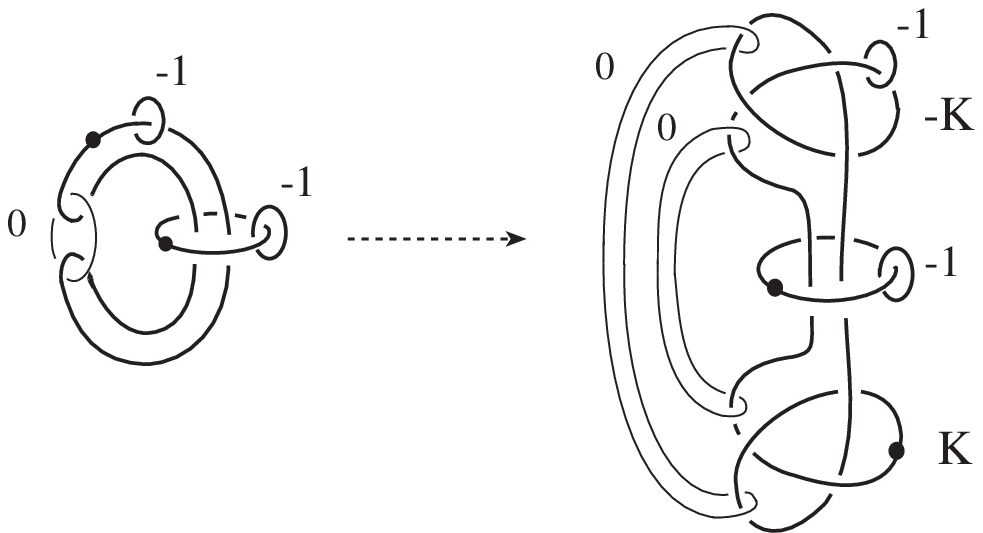}
\caption{}
\label{myfigure}
\end{figure}

Since the attaching circles of the other $2$-handles of $X$ could tangle to the boundary of $ T^{2}\times
D^2 $, it is important to indicate where the various linking circles of the boundary are thrown to by the
diffeomorphism of Figure 3. This is indicated in Figure 6. 

Recall that since  $3$- and $4$- handles of four manifolds are attached in the canonical way,
to describe a $4$-manifold it suffices to describe its $1$- and $2$- handle structure. So,
in order to visualize 
$ (S^{3} -K)\times S^1 $, which is obtained by by identifying
the two ends of $ (S^{3}- K)\times I$, it suffices to
visualize  
$ (B^{3} -K_{0})\times I $ with its ends identified, where
$K_{0}\subset B^{3}$ is a properly imbedded arc with the knot $K$ tied on it (the rest is a $3$-handle). The second
picture of  Figure 4 gives the handlebody picture of $ (B^{3} -K_{0} ) \times I$.

\begin{figure}[htb]
\includegraphics{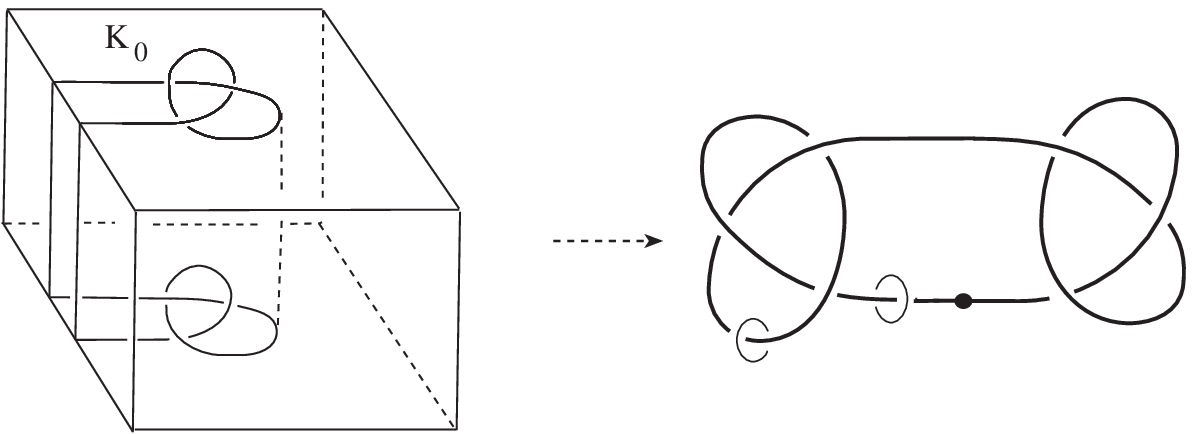}
\caption{}
\label{myfigure}
\end{figure}
 
Identifying the ends of $ (B^{3} -K_{0})\times I $
(up to 3-handles) corresponds to attaching a new 1-handle, and 2-handles, where the new $2$-handles 
are attached along the 1-handles of the two boundary components of $ (B^{3} -K_{0} ) \times I $ as indicated in
Figure 5 (more specifically the $2$-handles are attached along the loops connecting the {\it core circles} of the
knot complements).

\begin{figure}[htb]
\includegraphics{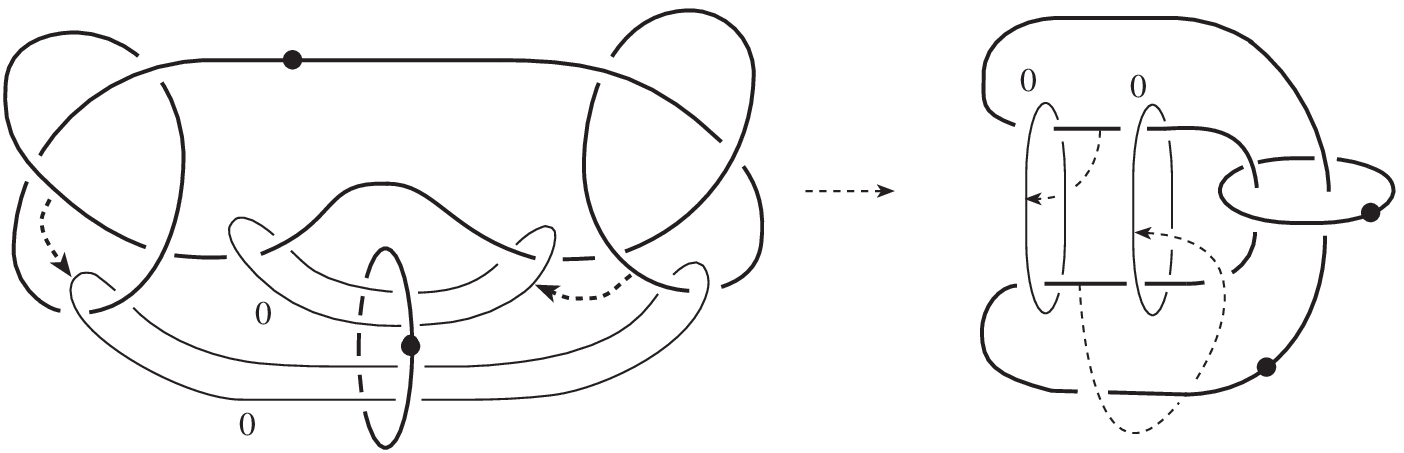}
\caption{}
\label{myfigure}
\end{figure}

To see the diffeomorphism of Figure 3 (i.e. to see that the boundary of 
the first picture in Figure 5
is standard), we simply remove the dot on the ``slice" $1$-handle  
(i.e. turn it to a $2$-handle) and
slide it over the two $2$-handles (as indicated by the arrows) in the first picture of 
Figure 5. This gives  the second picture of Figure
5. After sliding $2$-handles over each other of second picture of Figure 5, 
and cancelling the
resulting $ S^{2}\times D^{2} $ with the $3$-handle we obtain $ T^{2}\times D^{2} $. 
Also, to see the
inverse boundary diffeomorphism from $ T^{2}\times D^{2} $ to the 
first picture in Figure 5,  we remove the dot from
the $1$-handle of the second picture of Figure 5 and perform the  handle
slides indicated by the arrows. 

\begin{figure}[htb]
\includegraphics{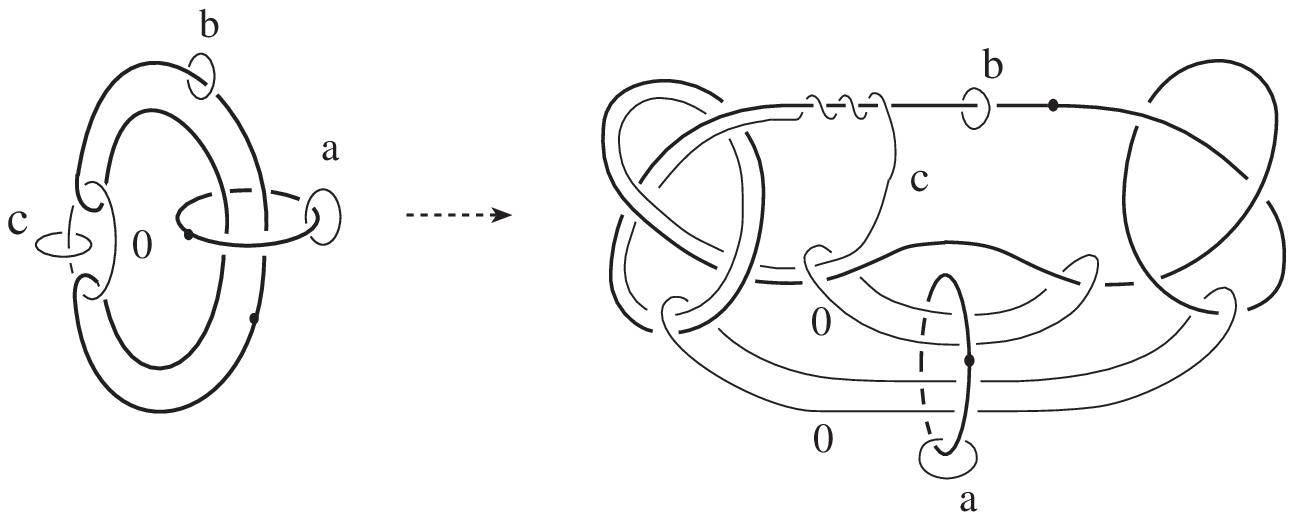}
\caption{}
\label{myfigure}
\end{figure}

Now putting these together in Figure 6 we see where the boundary diffeomorphism takes
various linking circles of $\partial (T^{2}\times D^{2})$. In particular the linking
circle $c$ of the $2$-handle is thrown to the loop which corresponds the zero push-off
of $K$  in $ K\# (-K) $.

Figure 7 is the same as the second picture of Figure 6 except that
the slice disk complement, which
$K \# (-K)$ bounds, is drawn more concretely. Also note that, though our discussion is for general
$K$, for the sake of concreteness, we have drawn our figures by taking $K$ to be the trefoil knot.

\begin{figure}[htb]
\includegraphics{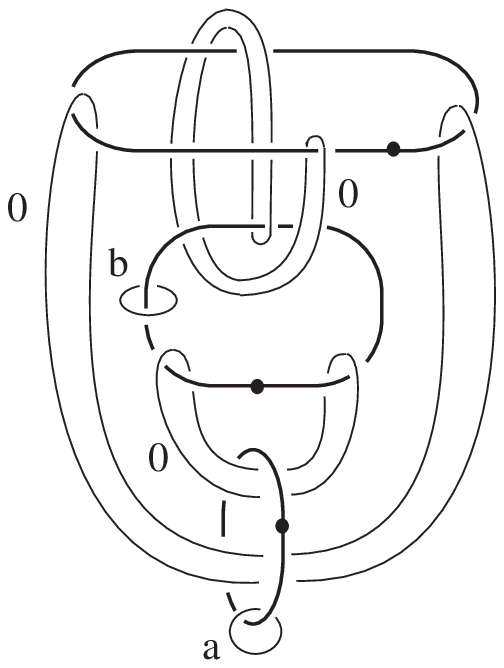}
\caption{}
\label{myfigure}
\end{figure}

\section{Applications}

In \cite {FS} Fintushel and Stern conjectured that if $X$ is the Kummer surface $K3$, 
then the
association 
$ K\to X_{K} $ gives an injective map from the isotopy classes of knots $K$  
in $ S^{3} $ to the
set of diffeomorphism classes of  smooth structures on $X$. 
The following provides a counterexample to
this conjecture:


\begin{thm}
 $X_{K} =X_{-K}$
\end{thm}

\begin{proof} 
 There is an obvious self-diffeomorphism of the second picture in Figure 3 
(i.e.  $(S^{3} -K)\times S^{1}$)
exchanging roles of
$K$ and $-K$ ; i.e. the diffeomorphism induced by $180 ^{o}$ rotation of ${\bf R}^{3}$ around the $y$-axis. It
is easily check that this diffeomorphism extends to the interior of $(S^{3}-K)\times S^{1}$, implying the
desired result.
\end{proof}

The following says that all smooth manifolds $X_{K}$
obtained from $X$ by using from different knots $K$ become standard after single stabilization. This result
was independently observed by Auckly \cite{Au}. 

\begin{thm}
 $X_{K} \;\# \;(S^{2}\times S^{2} ) =X\# (S^{2}\times S^{2} )$
\end{thm}

\begin{proof} 
  $X_{K} \;\# \;(S^{2}\times S^{2} ) $ is obtained by surgering any homotopically trivial loop (with the correct
framing). We choose to surger $X_{K}$ along the trivially linking circle of its slice $1$-handle (the knot  $ K\#
(-K) $ with a dot). This corresponds to turning the slice $1$-handle to a $0$-framed  $2$-handle (i.e. replace the dot
with $0$ framing), hence we are free to isotop the attaching circle of 
this $2$-handle  to the standard position as
indicated in Figure 5. In particular, 
this makes the boundary diffeomorphism between the two handlebodies of Figure 5
extend to a
$4$-manifold diffeomorphism. So, Surgered
$X_{K}$ is diffeomorphic to the surgered $X$ which is $X\# (S^{2}\times S^{2} )$. 

Note
that though we indicated the argument for the trefoil knot $K$ in our pictures, the
same applies for a general $K$ (i.e. in Figure 5 the knot $K \# (-K)$ unknots in the
presence of the
$2$-handles)
\end{proof}

Notice that $X_{K}$ can be viewed as
$ X_{K}= X_{f}= (X-T^2\times D^2)
\smile_{\varphi} (S^{3}\times S^{1}-U) $,
where $U$ is an open  tubular neighborhood of an imbedded torus 
$f : S^{1}\times S^{1} \to S^{3}\times S^{1} $, with
$Image (f )=K\times S^{1} $. The map $f$ is induced from the obvious 
imbedding $ K\times I
\to S^{3}\times I$ by identifying the ends. More generally to any 
concordance  $s$ from $K$ to itself, we can
associate an imbedding of a torus 
$f_{s}: S^{1}\times S^{1} \hookrightarrow S^{3}\times S^{1}$, hence getting map
$$ {\it C}(K) \longrightarrow \left\{ \begin {array}{c}
 \mbox {Diffeomorphism classes of} \\ \mbox {smooth structures on}\; X \end
 {array} \right\} $$

\noindent defined by $ s\to X_{s} $, where  $ {\it C}(K)= \{ \;  s: S^{1} \times I \hookrightarrow
B^{3}\times I\; |
\; s|_{S^{1}\times {0}}=s|_{S^{1}\times {1}} =K \;\} $.
It is an interesting question that  how the diffeomorphism class of
$X_{s}$ depends on the concordance class $s$  of $K$? The following says that the above map is not
injective.

\begin{thm}
If $K$ is the trefoil knot, there is $ s\in {\it C}(K\# (-K))$ 
such that  $X_{s} =X$
\end{thm}

\begin{proof} Let $s$ be the concordance of $ K\#(-K) $ to itself, given
by connected summing the two obvious slice discs which two copies of $ K\#(-K) $ bound
in
$B^{4}$ as in the second picture of Figure 8. 
\begin{figure}[htb]
\includegraphics{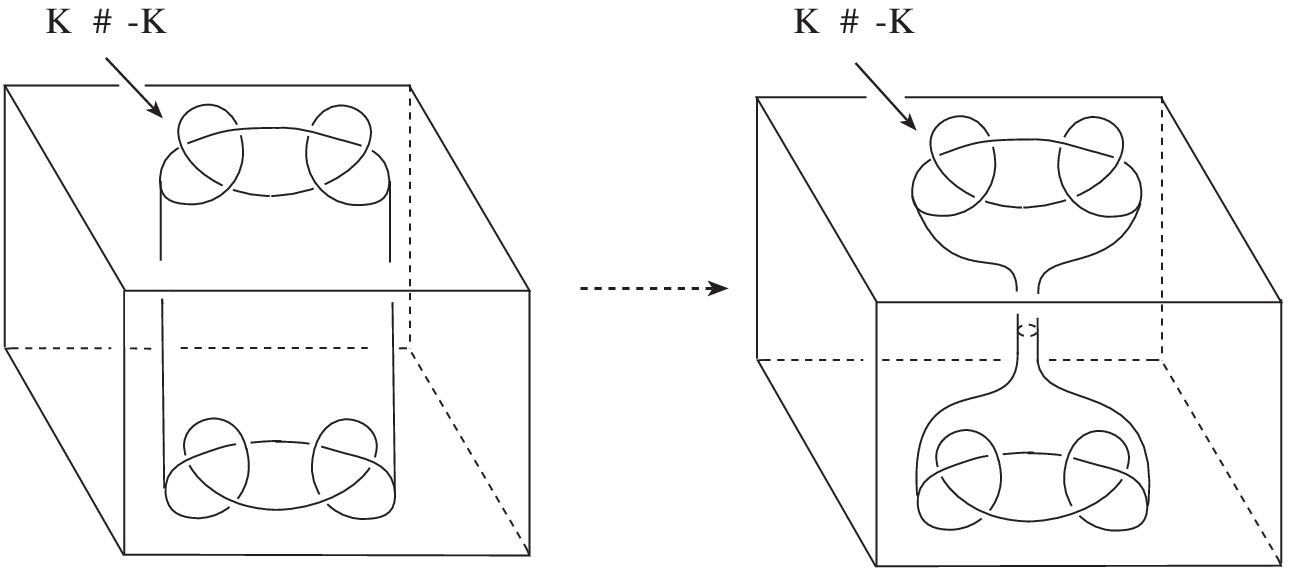}
\caption{}
\label{myfigure}
\end{figure}
Now if we use the product concordance $ s_{0}\; $ from $ K\#(-K) $ to itself,  i.e.
the first picture of Figure 8, our algorithm says that changing the  cusp
neighborhood by $ (S^{3} - K\#(-K))\times S^{1} $ is given by the handlebody of Figure 9, which is the same as
Figure 10 (where the slice $1$-handle is drawn as a usual handlebody). Whereas if we
use the concordance $s$, described above, we get Figure 11. By an isotopy we see that
Figure 11 is diffeomorphic  to Figure 12 which is diffeomorphic to Figure 13, and Figure 13 is isotopic to
Figure 14. By handle slides indicated  in Figures 14 and 15 we obtain the second picture of Figure 15. By
cancelling a 1-and 2- handle pairs we get the first  picture of Figure 16. Then by a
2-handle slide, and cancelling an unknotted 0-framed 2-handle by a 3-handle, 
we obtain the last picture of 16 which is
the cusp $C$. So we proved
$ C_{s}=C $, but since $C\subset X$ and every self
diffeomorphism of $ \partial C $ extends to $ C $ we conclude $X_{s}=X$  \end{proof}

\begin{rem}
 This theorem says that taking different elements $ s\in {\it C}(K) $ can result changing the  smooth structure
of
$X_{s}$. For example, if take any c-imbedded torus in a smooth manifold $ X $ with
$ SW_{X}\neq 0 $, and $K$ is the trefoil knot, and if $s_{0}\in {\it C}(K\# (-K))$
is the product concordance, then  by Theorem 1.1 $$ SW_{X_{s_{0}}} =
SW_{X_{K\#(-K)}}=SW_{X} . \Delta_{K\#(-K)} \neq  SW_{X}$$ hence $X_{s_{0}}\neq X$.
But  on the other hand Theorem 2.3 says that there is
$ s\in {\it C}(K\# (-K))$ with $ X_{s} =X $, so $X_{s_{0}} \neq X_{s}$. In particular, this shows that the
concordances $s_{0}$ and $s$ are different. This gives a hope the that hard to distinguish knot concordances 
might be distinguished by the Seiberg-Witten invariants of the associated manifolds $X_{s}$. 

\end{rem}

\begin{rem} Let $s \in C(K)$, and  $f_{s}: S^{1}\times S^{1} \hookrightarrow S^{3}\times S^{1}$ be the
 corresponding imbedding 
One can ask whether Theorem 1.1 generalizes to $ SW_{X_{s}} = SW_{X} .
\Delta_{s}(t) $?,  where $ \Delta_{s}(t) $ is the Alexander polynomial associated to this imbedding.
\end{rem}

\subsection {A twisted version of $X_{K}$}

Another version of the operation $X\to X_{K}$ that was previously introduced in \cite{CS},
which, in a sense, is the square root of this operation:
 Let $K$ is an invertible knot, i.e.   an orientation preserving involution $\tau: {\bf R}^{3}\to {\bf R}^{3}$
(e.g. $180^{0}$ rotation)  restricts to $K$ as
an involution with two fixed points, and let  $N$ be the open tubular neighborhood of
$K$.  Then we can form the following
$S^{3}-N$ bundle over $S^{1}$:
$$ (S^{3}- N)\tilde{\times} S^{1} = (S^{3}- N)\times [0,1]/ (x,0)\sim (\tau (x),1) $$
Define a $D^{2}$-bundle over the Klein bottle $ {\bf K} ^{2}$ by  
$C^{*} =S^{1}\times D^{2} \times [0,1] / (x,0)\sim (\tau (x),1)
$. Then $(S^{3}- N)\tilde{\times} S^{1}$ and $C^{*}$ have the same boundaries, and so if
$X$ is a smooth $4$-manifold with  $C^{*}\subset X$, we can construct
$$X_{K}^{*} = (X- C^{*}) \smile_{\varphi} (S^{3}-N)\tilde{\times} S^1$$ 

\noindent where $\varphi: \partial C^{*}\to \partial (K\times S^{1})\tilde{\times}S^{1} $ is a diffeomorphism with $
\; \varphi (p\times \partial D^2)=K\times p $. The operation 
$ X\to X^{*}_{K} $, is a certain generalization of
the  Fintushel and Stern operation $ X\to X_{K} $ done using a `Klein bottle'
instead of a torus. This operation was previously discussed in
\cite{CS}.  By using the previous arguments one can see see that the handlebody 
picture of the operation $ X\to X_{K}^{*} $ is given by  Figure 17. The first picture of Figure 17 is
$ C^{*} $ and the second is $ (S^{3}- N)\tilde{\times} S^{1} $. The rest of $
X_{K}^{*} $ is obtained by simply by drawing the images of the additional handles
under the diffeomorphism  $ \varphi : \partial  C^{*} \to \partial (S^{3}-
N)\tilde{\times} S^{1}$. For convenience, in Figure 17  the images of 
the linking circles $a,b,c$ under $\varphi$ are indicated. 

Now, call an imbedded Klein bottle
${\bf K}^{2}\subset X $  {\it c-imbedded Klein bottle}, if 
$${\bf K}^{2}\subset C^{*}\subset U\subset  X$$ 
where U is either one of the manifolds of Figure 18,
and $\pi_{1}(U)\to \pi_{1}(X)$ injects (notice $\pi_{1}(U)={\bf Z}_{2}$). Then it is
easy to see that the obvious $2$-fold cover $ \tilde{X}\to X $ contains a cusp $C$ (so it contains a {\it
c-imbedded }  $ T^{2} $), and the operation $ X\to X_{K}^{*} $  lifts to the
usual Fintushel-Stern knot surgery operation  $ \tilde{X}\to \tilde{X}_{K} $ (done
using this
$T^{2}$). Hence if  $ SW_{\tilde{X}}\neq 0 $ and $ \Delta_{K}(t) \neq 0
$, the operation $ X\to X^{*}_{K} $ changes the smooth structure of $ X $, i.e. $
X\neq X^{*}_{K} $. 
For example, $X$ can be a manifold with boundary, which is $2$-fold covered
 by a Stein manifold $ \tilde{X} $ (so $ \tilde{X} $ compactifies into a closed symplectic manifold
which Theorem 1.1 applies). It is easy to check that the first manifold of Figure 18 is such an example.

\begin{figure}[htb]
\includegraphics{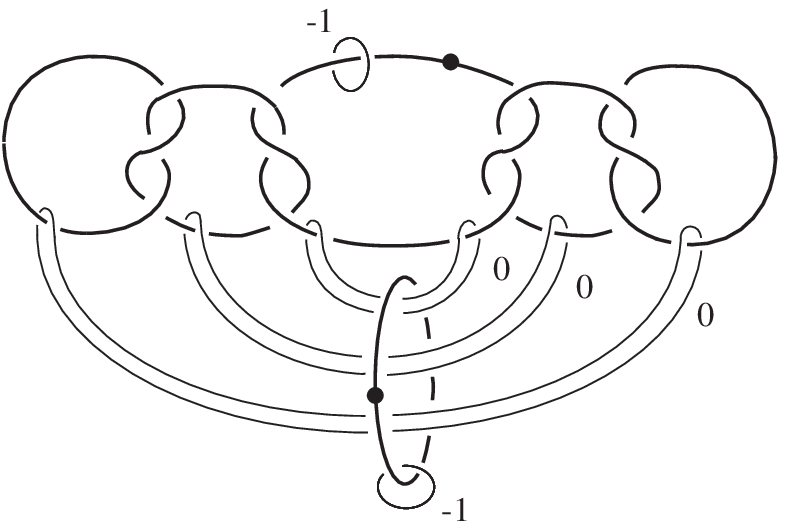}
\caption{}
\label{myfigure}
\end{figure}

\begin{figure}[htb]
\includegraphics{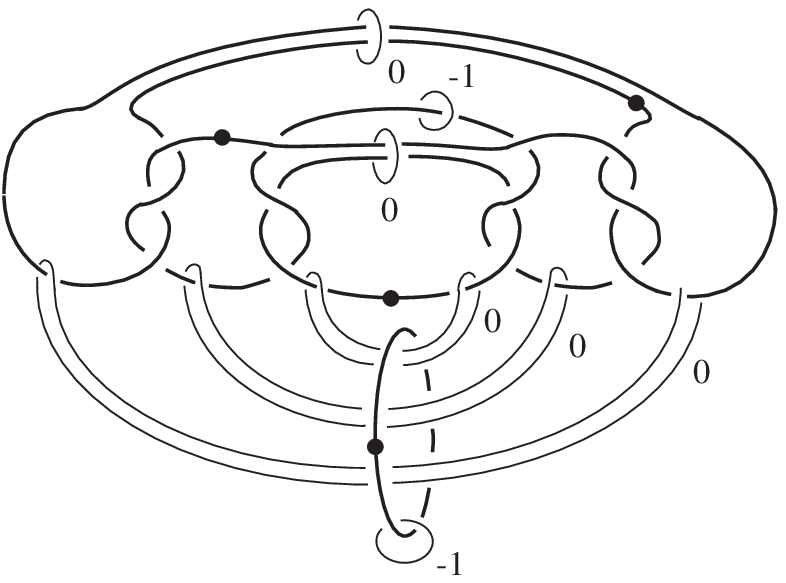}
\caption{}
\label{myfigure}
\end{figure}

\begin{figure}[htb]
\includegraphics{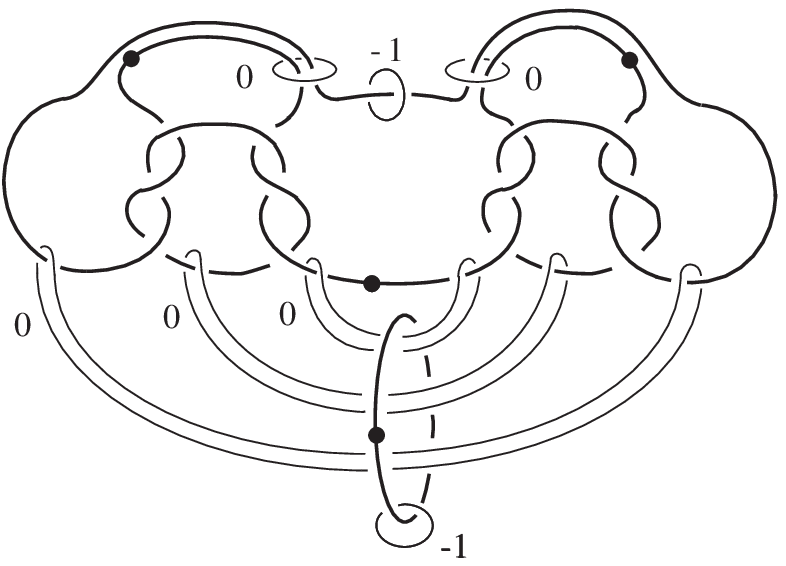}
\caption{}
\label{myfigure}
\end{figure}

\begin{figure}[htb]
\includegraphics{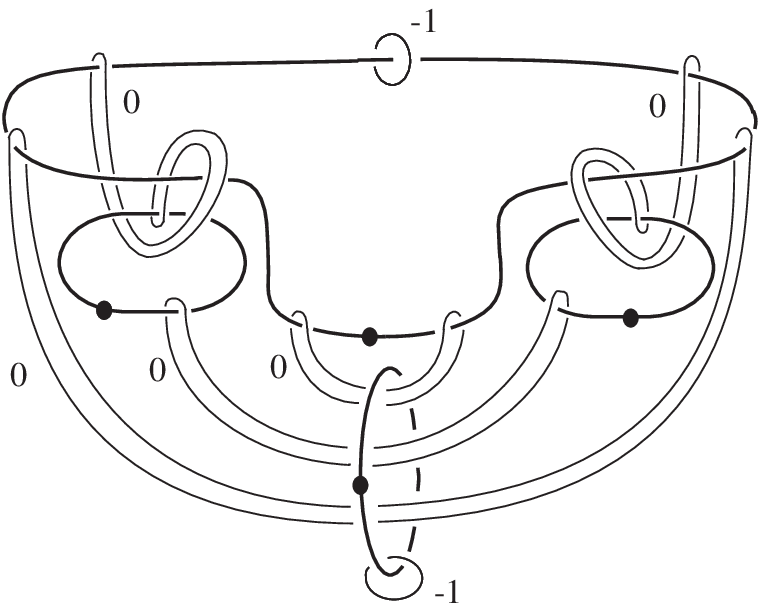}
\caption{}
\label{myfigure}
\end{figure}

\begin{figure}[htb]
\includegraphics{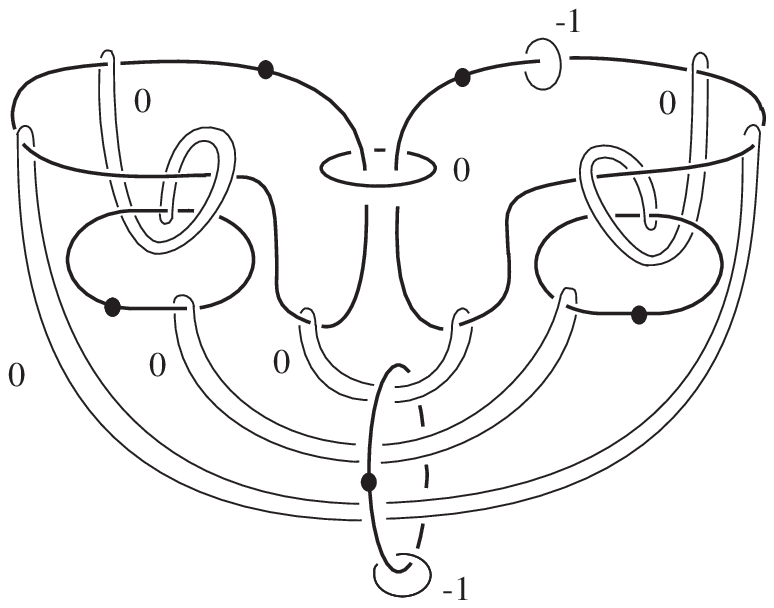}
\caption{}
\label{myfigure}
\end{figure}

\begin{figure}[htb]
\includegraphics{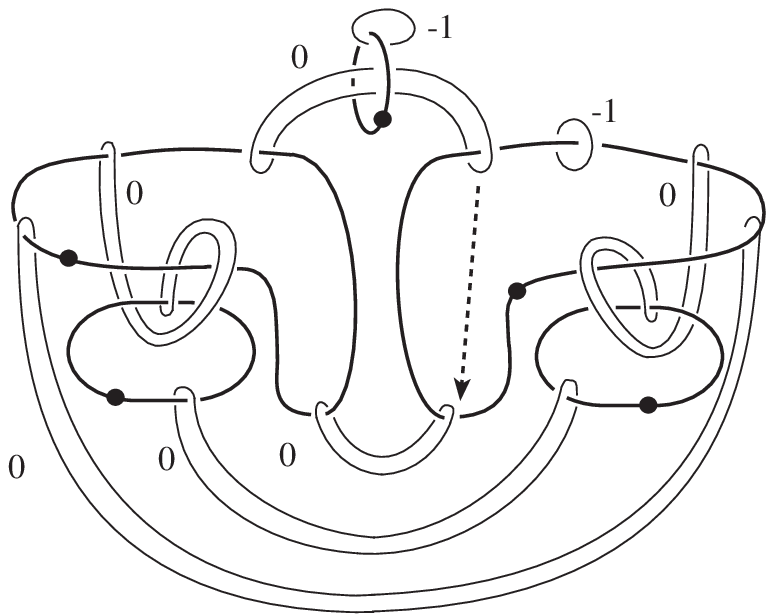}
\caption{}
\label{myfigure}
\end{figure}

\begin{figure}[htb]
\includegraphics{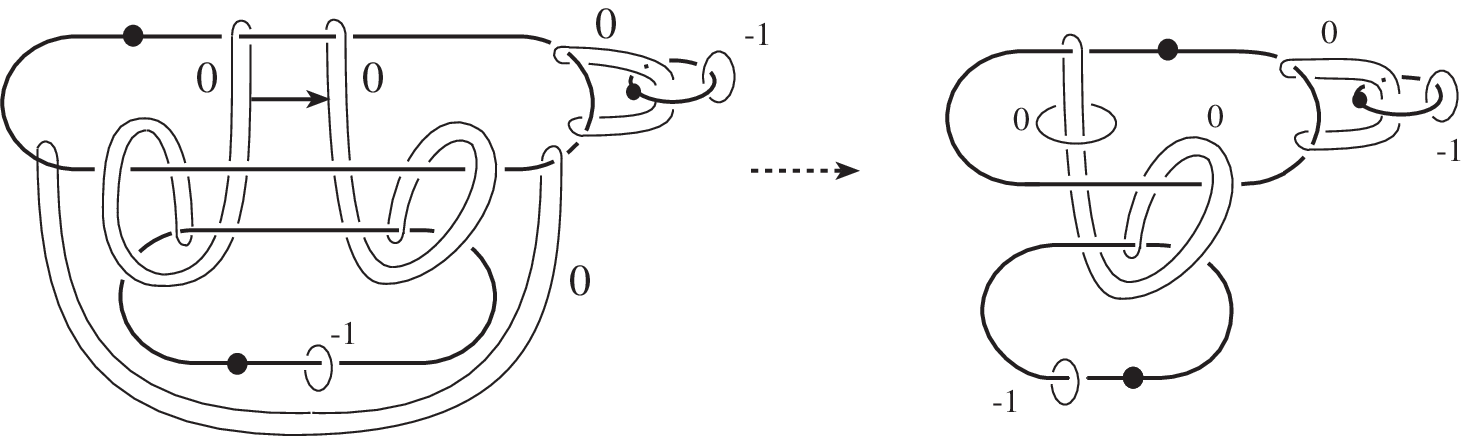}
\caption{}
\label{myfigure}
\end{figure}

\begin{figure}[htb]
\includegraphics{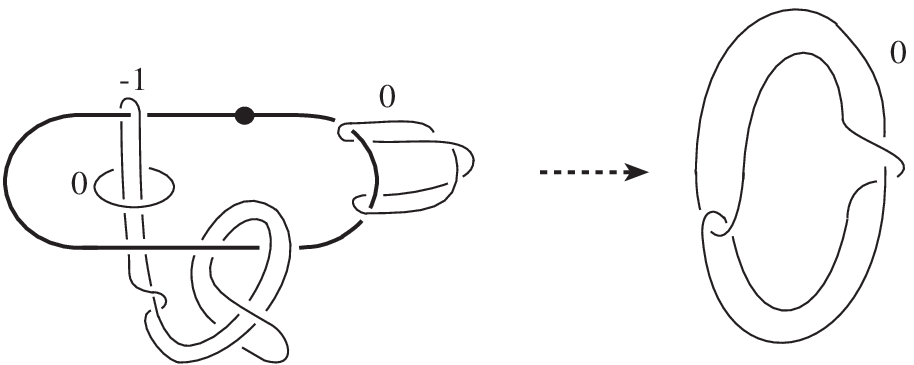}
\caption{}
\label{myfigure}
\end{figure}

\begin{figure}[htb]
\includegraphics{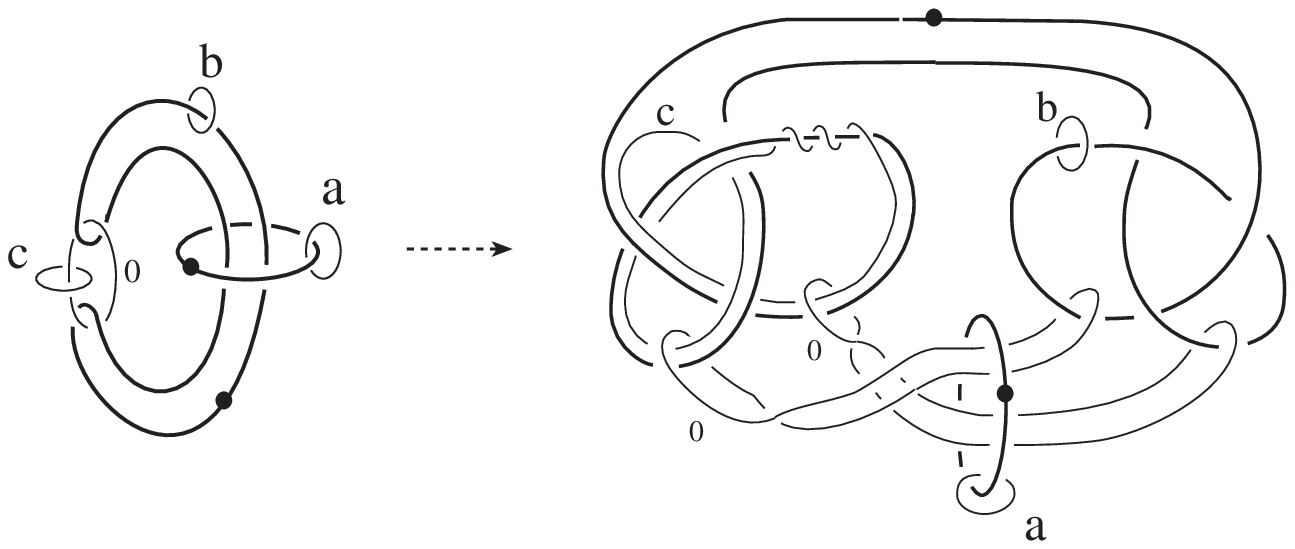}
\caption{}
\label{myfigure}
\end{figure}

\begin{figure}[htb]
\includegraphics{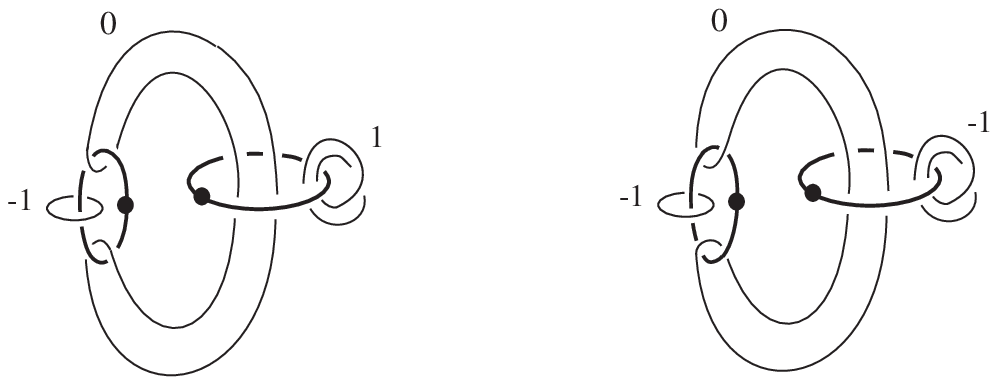}
\caption{U}
\label{myfigure}
\end{figure}

\end{document}